\documentclass{amsart}

\usepackage[T1]{fontenc}
\usepackage[british]{babel}
\usepackage{a4}
\usepackage{amssymb}
\usepackage{diagrams}

\makeatletter

\DeclareMathAlphabet\mathbfit{OML}{cmm}{b}{it}

\def\comment#1{\M@margin{#1}}

\def\ie{i.\,e.}

\def\cf{cf.}

\def \N{{\bf N}}
\def \Z{{\bf Z}}

\def\degree#1{\lvert#1\rvert}

\DeclareMathOperator{\id}{id}

\DeclareMathOperator{\Hom}{Hom}
\DeclareMathOperator{\Tor}{Tor}
\DeclareMathOperator{\Sym}{Sym}

\def\g{\mathfrak g}

\def\M@warning#1{\typeout{MATH Warning: #1.}}

\def\M@margin#1{\if@printlabels\setbox0=\vbox to\z@
  {\vss\hbox to\z@{\hskip\hsize\hskip\labelskip\footnotesize #1\hss}}%
  \dp0=\z@\ifvmode\box0\else\vadjust{\box0}\fi\fi}

\def \M@restore{\catcode `\^=7 \catcode`\_=8 }

\def \newterm{\@ifnextchar[{\y@newterm}{\x@newterm}}
\def \x@newterm #1{\y@newterm[#1]{#1}}
\def \y@newterm[#1]#2{
  \textbf{\boldmath{#2}}}

\def \newsymbol #1{\M@margin{$\mapsto #1$}#1}

\def\printsymbols{{\def\indexname{Symbols}%
  \begin{theindex}%
  \normalsize
  \makeatletter
  \@input{\jobname.sym}%
  \makeatother
  \end{theindex}}}

\expandafter \ifx \csname diagram\endcsname \relax \else
\let \o@diagram=\diagram
\@namedef{diagram*}{\M@restore \o@diagram}
\@namedef{enddiagram*}{\enddiagram}
\def \diagram{\M@restore \refstepcounter{equation}%
        \o@diagram[eqno=\@eqnnum,moreoptions]}
\fi

\def\indexhat#1{\M@margin{$\to$~#1}\index{#1}#1}
\let\o@sp\sp
\def\sp{\ifmmode\let\M@temp\o@sp\else\let\M@temp\indexhat\fi\M@temp}
\catcode`\^=\active
\def^{\protect\sp}

\newif\if@printlabels
\let\printlabels=\@printlabelstrue
\let\printnolabels=\@printlabelsfalse
\newdimen\labelskip
\printnolabels
\let\o@label=\label
\newtoks\M@everypar
\def\M@label#1{\M@margin{\tt #1}\o@label{#1}%
  {\def\@currentlabel{\@currentthm}\o@label{lthm@#1}}%
  \@ifundefined{full@label}{}{\def\@currentlabel{\full@label}%
  \o@label{fthm@#1}}}
\def\label#1{\@ifundefined{@currentthm}{\o@label{#1}}%
    {\ifvmode\M@everypar=\everypar\everypar={\the\M@everypar\M@label{#1}}%
     \else\M@label{#1}\fi}}
\printnolabels

\let\o@enumerate\enumerate
\def\enumerate{\edef\full@label{%
  \@ifundefined{full@label}{\@currentlabel}{\full@label\,(\@currentlabel)}}%
  \o@enumerate}

\@ifpackageloaded{amsthm}{%
  \let\o@ynthm=\@ynthm  
  \def\@ynthm#1[#2]#3{\global\@namedef{dthm@#1}{#3}\o@ynthm{#1}[#2]{#3}}%
  }{%
  \let\o@othm=\@othm  
  \let\o@nthm=\@nthm
  \def\@othm#1[#2]#3{\global\@namedef{dthm@#1}{#3}\o@othm{#1}[{#2}]{#3}}%
  \def\@nthm#1#2{\global\@namedef{dthm@#1}{#2}\o@nthm{#1}{#2}}%
  }%
\let \o@thm=\@thm
\def \@thm{\let\@currentthm=\@currenvir\o@thm}
\def\@thmwarning#1{\M@warning{(no)theoremref: Reference `#1'
        on page \thepage \space not to a theorem}}
\def \@thmprintnum#1{\@ifundefined{r@fthm@#1}%
  {\ref{#1}}{\ref{fthm@#1}\,(\ref{#1})}}
\def\@thmcheck#1[#2]#3{\@ifundefined{r@#3}%
        {{\bf ??}\M@warning{(no)theoremref: Reference `#3' on page
        \thepage \space undefined}}%
        {\@ifundefined{r@lthm@#3}{\@thmwarning{#3}}%
        {\@ifundefined{dthm@#2}{\M@warning{(no)theoremref: Theorem
        environment `#2' on page \thepage \space undefined}}%
        {\edef \@tempa{\@nameuse{r@lthm@#3}}%
        \edef \@tempb{\expandafter \@car \@tempa \@nil}%
        \def \@tempa{#2}%
        \ifx \@tempa \@tempb \else
        \M@warning{(no)theoremref: Reference `#3' on page \thepage \space not
        to a `#2'}\fi}%
        #1\@thmprintnum{#3}\fi}}}
\def\@thmprintall#1{\@ifundefined{r@#1}%
        {{\bf ??}\M@warning{Reference `#1' on page \thepage \space undefined}}%
        {\x@thmprintall{#1}}}
\def \x@thmprintall#1{\@ifundefined{r@lthm@#1}%
        {\@thmwarning{#1}}{\y@thmprintall{#1}}}
\def \y@thmprintall#1{\edef \M@temp{\@nameuse{r@lthm@#1}}%
        \@nameuse{dthm@\expandafter \@car \M@temp \@nil}~\@thmprintnum{#1}}
\def\theoremref{\@ifnextchar[{\@thmcheck \iftrue}{\@thmprintall}}
\def\notheoremref{\@thmcheck{\iffalse}}

\def\proofsectionref#1{{\bf ??}}

\def\M@gobblefour#1#2#3#4{}
\def\proofsection[#1]#2{%
  \section[Proof of #1]{Proof of \protect\theoremref{#2}}\label{p@#2}}

\def\smashsubstack#1{\hbox to5ex{\hss$\substack{#1}$\hss}}
\def\sumsss#1{\sum_{\smashsubstack{#1}}}

\def\cleardoublepage{\clearpage\if@twoside \ifodd\c@page\else
    \hbox{}\thispagestyle{empty}\newpage\if@twocolumn\hbox{}\newpage\fi\fi\fi}

\newtheorem{theorem}{Theorem}[section]
\newtheorem{lemma}[theorem]{Lemma}
\newtheorem{proposition}[theorem]{Proposition}

\newenvironment{acknowledgements}{\begin{trivlist}\item[]{\em Acknowledgements.}\hskip0.5em}{\end{trivlist}}

\makeatother

\let\epsilon\varepsilon

\newtheorem{remark}[theorem]{Remark}

\numberwithin{equation}{section}
\allowdisplaybreaks[1]

\newarrow{Equal}=====

\def\and{\qquad\text{and}\qquad}

\let\alg\mathbf
\let\Sim\mathbfit

\def\oneptspace{{\rm pt}}

\def\Ll{\boldsymbol\Lambda}
\def\Sl{\mathbf S}

\def\Su{\Sl^*}
\def\Lu{\Ll^{\!*}}

\def\LMod{\hbox{\bf\em M\kern-0.05em od}}
\def\lMod#1{#1\hbox{-\bf\em M\kern-0.05em od}}

\def\rComod#1{\hbox{{\bf\em Comod}-}#1}
\def\wMod#1{#1\hbox{-$\mathcal M$\bf\em\kern-0.1em od}}
\def\lComod#1{#1\hbox{-\bf\em C\kern-0.05em omod}}
\def\wComod#1{\rComod{B\Ll}}

\let\shuffle\nabla
\def\AW{AW}

\def\cuponeproduct{cup$_1$~product}
\def\Bf{\psi}

\def\tcomp{\Psi}

\def\Space#1{#1\hbox{-{\bf\em Space}}}
\def\Spaceover#1{\hbox{{\bf\em Space}-}{#1}}

\def\timesunder#1{\mathbin{\mathchoice
  {\mathop\times\limits_{\mkern-5mu #1\mkern-5mu}}%
  {\times_{#1}}{\times_{#1}}{\times_{#1}}}}

\def\otimesunder #1{\mathbin{\mathchoice
  {\mathop\otimes\limits_{\mkern-20mu #1\mkern-20mu}}%
  {\otimes_{#1}}{\otimes_{#1}}{\otimes_{#1}}}}

\def\g{\mathfrak g}
\def\Ainfty{A_\infty}
\let\shmto\Rightarrow

\title{Koszul duality and equivariant cohomology}
\author{Matthias Franz}

\address{Section de Math\'ematiques, Universit\'e de Gen\`eve,
  Case postal 240, 1211 Gen\`eve 24, Switzerland}
\email{matthias.franz@math.unige.ch}


\begin{document}

\begin{abstract}
  Let $G$ be a topological group such that its homology~$H(G)$
  with coefficients in a principal ideal domain~$R$ is an exterior algebra,
  generated in odd degrees.
  We show that the singular cochain functor carries
  the duality between $G$-spaces and spaces over~$BG$
  to the Koszul duality between modules up to homotopy
  over $H(G)$~and~$H^*(BG)$.
  This gives in particular a Cartan-type model
  for the equivariant cohomology of a $G$-space
  with coefficients in~$R$.
  As another corollary, we obtain a multiplicative
  quasi-isomorphism~$C^*(BG)\to H^*(BG)$.

  A key step in the proof is to show that a differential Hopf algebra
  is formal in the category of $\Ainfty$~algebras
  provided that it is free over~$R$ and its homology an exterior algebra.
\end{abstract}

\maketitle

\section{Introduction}


Let $G$ be a topological group.
A space over the classifying space~$BG$ is a map~$Y\to BG$.
There are canonical ways to pass from left $G$-spaces
to spaces over~$BG$ and back:
The Borel construction~$\Sim t X =EG\timesunder G X$ is a functor
\[
  \Sim t\colon\Space G\to\Spaceover{BG},
\]
and pulling back the universal right $G$-bundle~$EG\to BG$ along~$Y\to BT$
and passing to a left action gives a functor in the other direction,
\[
  \Sim h\colon\Spaceover{BG}\to\Space G.
\]
These functors are essentially inverse to each other in the sense that
$\Sim h\Sim t X$~and~$\Sim t\Sim h Y$ are homotopy-equivalent
in the category of spaces to $X$~and~$Y$, respectively,
\cf~\cite{DrorDwyerKan:80}. 

We show that
the normalised singular cochain functor~$C^*$ transforms this duality,
up to quasi-isomorphism,
to the Koszul duality between modules up to homotopy
over the homology~$\Ll=H(G)$ and the cohomology~$\Su=H^*(BG)$.
The only assumptions are that coefficients
are in a principal ideal domain~$R$ and that $H(G)$ is an exterior algebra
on finitely many generators of odd degrees
or, equivalently, $H^*(BG)$ a symmetric algebra
on finitely many generators of even degrees.

A priori, the isomorphism~$H(G)\cong\bigwedge(x_1,\ldots,x_r)$
must be one of Hopf algebras\footnote{%
Note that $H(G)$ has a well-defined diagonal because it is free over~$R$.}
with primitive generators~$x_i$.
But the Samelson--Leray theorem asserts that in our situation
any isomorphism of algebras (or coalgebras)
can be replaced by one which is Hopf.
In characteristic~$0$ it suffices by Hopf's theorem to check
that $G$~is connected and $H(G)$ free of finite rank over~$R$.
In particular, the condition is satisfied
for $U(n)$,~$SU(n)$ and~$Sp(n)$ and arbitrary~$R$,
and for an arbitrary compact connected Lie group
if the order of the Weyl group is invertible in~$R$.

Previously,
Goresky, Kottwitz and MacPherson~\cite{GoreskyKottwitzMacPherson:98}
related the  topological and Koszul duality for real coefficients,
compact connected Lie groups and subanalytic spaces.
Since they used differential forms instead of singular cochains,
their complexes were strict modules and not only ones up to homotopy.
For tori, the generalisation to the present setting was achieved
by the author in~\cite{Franz:03}. There it was also shown how to extend
the Koszul duality functors
$$
  \alg t\colon\lMod\Ll\to\lMod\Su
$$
and
$$
  \alg h\colon\lMod\Su\to\lMod\Ll
$$
to modules up to homotopy
(or ``weak modules'', as we will call them).
For torus actions there is actually no need to consider weak $\Ll$-modules
because the $\Ll$-action on cohomology lifts to an honest action on cochains.
It is not clear how to do this
for non-commutative groups.\footnote{In~%
\cite[Sec.~12]{GoreskyKottwitzMacPherson:98}
it is claimed that such a lifting is possible
for any compact connected Lie group,
but the proof is wrong.
The mistake is that it is not possible in general
to find conjugation-invariant representatives of the generators~$x_i$,
see the footnote to Lemma~4.1 in~\cite{Franz:03}.}

Recall that a weak $\Ll$-module
is a differential module~$N$ over some differential algebra~$A$
together with a twisting cochain~$\Sl\to A$.
Its homology~$H(N)$ has a well-defined action of~$\Ll$.
Similarly, a weak $\Su$-module is an $A$-module~$M$
together with a twisting cochain~$\Lu\to A$,
and its homology is canonically an $\Su$-module.
Here $\Lu=H^*(G)$~and~$\Sl=H(BG)$ are the coalgebras dual
to~$\Ll$~and~$\Su$, respectively.
(Twisting cochains will be reviewed in Section~\ref{preliminaries}.)

\begin{proposition}\label{existence-twisting-cochains}
  There are twisting cochains~$v\colon\Sl\to C(G)$ and~$u\colon \Lu\to C^*(BG)$
  such that the $\Ll$-action on the homology of a $C(G)$-module,
  viewed as weak $\Ll$-module, is the canonical one over~$H(G)=\Ll$,
  and analogously for~$u$.
\end{proposition}

The cochains on a $G$-space are canonically a
$C(G)$-module and the cochains on a space over~$BG$ a
$C^*(BG)$-module.
Hence we may consider $C^*$ as a functor from $G$-spaces
to weak $\Ll$-modules, and from spaces over~$BG$ to weak $\Su$-modules.

We say that two functors to a category of complexes are quasi-isomorphic
if they are related by a zig-zag of natural transformations
which become isomorphisms after passing to homology.

\begin{theorem}\label{main-result}
  The functors $C^*\circ\Sim t$~and~$\alg t\circ C^*$
  from $G$-spaces to weak $\Su$-modules are quasi-isomorphic,
  as are the functors $C^*\circ\Sim h$~and~$\alg h\circ C^*$
  from spaces over~$BG$ to weak $\Ll$-modules.
\end{theorem}

Hence, the equivariant cohomology~$H_G^*(X)$ of a $G$-space~$X$
is naturally isomorphic, as $\Su$-module,
to the homology of the ``singular Cartan model''
\begin{subequations}\label{complex-Cartan}
\begin{equation}
   \alg t C^*(X)=\Su\otimes C^*(X) 
\end{equation}
with differential
\begin{equation}
  d(\sigma\otimes\gamma)
  =\sigma\otimes d\gamma
  +\sum_{i=1}^r\xi_i\sigma\otimes c_i\cdot\gamma
  +\sum_{i\le j}\xi_i\xi_j\sigma\otimes c_{i j}\cdot\gamma
  +\cdots,  
\end{equation}
\end{subequations}
where the $\xi_i$ are generators of the symmetric algebra~$\Su$
and the~$c_i\in C(G)$ representatives of the generators~$x_i\in\Ll$.
They are, like the higher order terms~$c_{i j}$ etc., encoded
in the twisting cochain~$v$. The sum, which runs over all non-constant
monomials of~$\Su$, is well-defined for degree reasons.

Similarly, the cohomology of the pull back of~$EG$ 
along~$Y\to BG$ is isomorphic to the homology of the
$\Ll$-module~$\alg h C^*(Y)=\Lu\otimes C^*(Y)$,
again with a twisted differential.
(See Section~\ref{Koszul-duality} for precise formulas for the differentials.)
That the complex~$\alg h C^*(Y)$ gives the right cohomology as $R$-module
is due to Gugenheim--May~\cite{GugenheimMay:74},
who also constructed a twisting cochain~$\Lu\to C^*(BG)$.
The correctness of the $\Ll$-action is new,
as is the singular Cartan model~\eqref{complex-Cartan}.

Along the way, we obtain the following result,
which was previously only known for tori,
and for other Eilenberg--Mac\,Lane spaces if~$R=\Z_2$
(Gugenheim--May~\cite[\S4]{GugenheimMay:74}):

\begin{proposition}\label{multiplicative-quasi-isomorphism}
  There exists a quasi-isomorphism of algebras~$C^*(BG)\to H^*(BG)$
  between the cochains and the cohomology of the simplicial construction
  of the classifying space of~$G$.
\end{proposition}

Any such map has an $\Ainfty$ map as homotopy inverse
(\cf~Lemma~\ref{shm-inverse}).
So we get as another corollary the well-known existence
of an $\Ainfty$ quasi-isomorphism $H^*(BG)\shmto C^*(BG)$.
The original proof (Stasheff--Halperin~\cite{StasheffHalperin:70})
uses the homotopy-com\-mu\-ta\-tiv\-ity of the cup product
and the fact that $H^*(BG)$ is free commutative.
Here it is based, like most of the paper, on the following result,
which is of independent interest and
should be considered as dual to the theorem of Stasheff and Halperin.

\begin{theorem}\label{A-infinity-formal}
  Let $A$ be a differential $\N$-graded Hopf algebra,
  free over~$R$ and such that
  its homology is an exterior algebra
  on finitely many generators of odd degrees.
  Then there are $\Ainfty$~quasi-isomorphisms
  $A\shmto H(A)$~and~$H(A)\shmto A$.
\end{theorem}

It is essentially in order to use Theorem~\ref{A-infinity-formal}
(and a similar argument in Section~\ref{end-other-proofs})
that we assume $R$ to be a principal ideal domain.
A look at the proofs will show that once
Proposition~\ref{existence-twisting-cochains},
Theorem~\ref{main-result} and Proposition~\ref{multiplicative-quasi-isomorphism}
are established for such an~$R$, they follow by extension of scalars
for any commutative $R$-algebra~$R'$ instead of~$R$.


Johannes Huebschmann has informed the author
that he has been aware of Theorem~\ref{A-infinity-formal} since the 1980's.
Instead of adapting arguments
from his habilitation thesis~\cite[Sec.~4.8]{Huebschmann:84},
we shall base the proof on an observation due to Stasheff~\cite{Stasheff:63b}.

\medbreak

The paper is organised as follows:
Notation and terminology is fixed in Section~\ref{preliminaries}.
Section~\ref{Koszul-duality} contains a review of Koszul duality
between modules up to homotopy over symmetric and exterior algebras.
Theorem~\ref{A-infinity-formal} is proved in Section~\ref{proof-formal}.
The proofs of the other results stated in the introduction appear
in Sections \ref{start-other-proofs}~to~\ref{end-other-proofs}.
In Section~\ref{equivariant-cohomology} we discuss equivariantly formal spaces
and in Section~\ref{other-models}
the relation between the singular Cartan model and other models,
in particular the classical Cartan model.
In an appendix we prove the versions of the theorems of Samelson--Leray
and Hopf mentioned above because they are not readily available
in the literature.

\begin{acknowledgements}
  The author thanks St\'ephane Guillermou, Johannes Huebsch\-mann,
  Tomasz Maszczyk and Andrzej Weber
  for stimulating discussions and comments.
\end{acknowledgements}

\section{Preliminaries}\label{preliminaries}

Throughout this paper, the letter~$R$ denotes a principal ideal domain.
All complexes are over~$R$.
Differentials always lower degree,
hence cochain complexes and cohomology are negatively graded.
All (co)algebras and (co)modules are graded and have differentials
(which might be trivial).
Let $A$~and~$B$ be complexes.
The dual~$f^*\in\Hom(B^*,A^*)$ of a map~$f\in\Hom(A,B)$ is defined by
$$
  f^*(\beta)(a)=(-1)^{\degree f\degree\beta}\,\beta(f(a)).
$$

Algebras will be associative and coalgebras coassociative,
and both have units and counits (augmentations).
Morphisms of (co)algebras preserve these structures.
We denote the augmentation ideal of an algebra~$A$ by~$\bar A$.
An $\N$-graded algebra~$A$ is called connected if $\bar A_0=0$,
and an $\N$-graded coalgebra~$C$ simply connected if $C_0=R$~and~$C_1=0$.
Hopf algebras are algebras which are also coalgebras
with a multiplicative diagonal, \cf~\cite[Def.~4.39]{McCleary:01}.
(Note that we do not require the existence of an antipode,
though there will always be one for our examples.)

Let $C$ be a coalgebra, $A$ an algebra and~$t\colon C\to A$ a twisting cochain.
For a right $C$-comodule~$M$ and a left $A$-module~$N$, we define the
twisted tensor product~$M\otimes_t N$ with differential
$$
  d_t=
  d\otimes1+1\otimes d+(1\otimes\mu_N)(1\otimes t\otimes 1)(\Delta_M\otimes1).
$$
Here $\Delta_M\colon M\to M\otimes C$ and $\mu_N\colon A\otimes N\to N$
denote the structure maps of $M$~and~$N$, respectively.
Readers unfamiliar with twisting cochains can take the fact
that $d$ is a well-defined differential (say, on~$C\otimes_t A$)
as the definition of a twisting cochain,
plus the normalisation conditions $t\iota_C=0$~and~$\epsilon_A t=0$, where
$\iota_C$ is the unit of~$C$ and $\epsilon_A$ the augmentation of~$A$.
Suppose that $C$~and~$A$ are $\N$-graded.
We will regularly use the fact that twisting cochains~$C\to A$
correspond bijectively to coalgebra maps~$C\to BA$
and to algebra maps~$\Omega C\to A$.
Here $BA$ denotes the normalised bar construction of~$A$ and $\Omega C$
the normalised cobar construction of~$C$.
In particular, the functors $\Omega$~and~$B$ are adjoint.
(See for instance \cite[Sec.~II]{HusemollerMooreStasheff:74}
for more about twisting cochains and the (co)bar construction.)

We agree that an exterior algebra is one on finitely many generators
of odd positive degrees. Let $A$ be an $\N$-graded algebra such that
$\Ll=H(A)=\bigwedge(x_1,\ldots,x_r)$ is an exterior algebra.
Then $H(B A)=H(B\Ll)=\Sl$ is a symmetric coalgebra on finitely many
cogenerators~$y_i$ of even degrees~$\degree{y_i}=\degree{x_i}+1$,
\cf~\cite[Thm.~7.30]{McCleary:01}.
(The converse is true as well.)
We assume that the~$y_i$ are chosen such that they
can be represented by the cycles~$[x_i]\in B\Ll$ and~$[c_i]\in B A$,
where the~$c_i\in A$ are any representatives of the generators~$x_i\in\Ll$.
We denote by~$x_\pi$, $\pi\subset\{1,\ldots,r\}$, the canonical $R$-basis
of~$\Ll$ generated by the~$x_i$, and the dual basis of~$\Lu$ by~$\xi_\pi$.
The $R$-basis of~$\Sl$ induced by the~$y_i$
is written as~$y_\alpha$, $\alpha\in\N^r$.
The dual~$\Su$ of~$\Sl$ is a symmetric algebra on generators~$\xi_i$
dual to the~$y_i$.

We work in the simplicial category.
We denote by~$C(X)$ the normalised chain complex of the simplicial set~$X$.
(If $X$ comes from a topological space,
then $C(X)$ is the complex of normalised singular chains.)
The (negatively graded) dual complex of normalised cochains
is denoted by~$C^*(X)$.
If $G$ is a connected (topological or simplicial) group,
then the inclusion of the simplicial subgroup consisting of the simplices
with all vertices at~$1\in G$ is a quasi-isomorphism.
We may therefore assume that $G$ has only one vertex.
Then $C(G)$ is a connected Hopf algebra
and $C(BG)$ a simply connected coalgebra.
In both cases, the diagonal is the Alexander--Whitney map,
and the Pontryagin product of~$C(G)$ is the composition of the
shuffle map~$C(G)\otimes C(G)\to C(G\times G)$
with the map~$C(G\times G)\to C(G)$ induced by the multiplication of~$G$.
Analogously, $C(X)$ is a left $C(G)$-module if $X$ is a left $G$-space.
The left $C(G)$-action on cochains is defined by
\begin{equation}\label{action-cochains}
  (a\cdot\gamma)(c)=(-1)^{\degree a\degree\gamma}\,\gamma(\lambda_*(a)\cdot c)
\end{equation}
where~$\lambda\colon G\to G$ denotes the group inversion.
If $p\colon Y\to BG$ is a space over~$BG$, then $C^*(Y)$
is a left $C^*(BG)$-module by $\beta\cdot\gamma=p^*(\beta)\cup\gamma$.

\section{Koszul duality}\label{Koszul-duality}

Koszul duality is most elegantly expressed as a duality
between $\Ll$-modules and comodules over the symmetric coalgebra~$\Sl$
dual to~$\Su$, see~\cite[Sec.~2]{Franz:03}.
It hinges on the fact that the Koszul complex~$\Sl\otimes_w\Ll$ is acyclic,
where $w\colon\Sl\to\Ll$ is the canonical twisting cochain
which sends each~$y_i$ to~$x_i$ and annihilates all other~$y_\alpha$.
In this paper, though, we adopt a cohomological viewpoint.
This makes definitions look rather ad hoc,
but it is better suited to our discussion of equivariant cohomology
in Section~\ref{equivariant-cohomology}.

We denote the categories of bounded above weak modules over $\Ll$~and~$\Su$
by $\wMod\Ll$~and~$\wMod\Su$, respectively.
(Recall that we grade cochain complexes negatively.)
Note that any (strict) module over $\Ll$~or~$\Su$ is also a weak module
because of the canonical twisting cochain~$w$
and its dual~$w^*\colon\Lu\to\Su$.
The homology of a weak $\Ll$-module~$(N,v)$ is a $\Ll$-module by
setting $x_i\cdot[n]=[v(y_i)\cdot n]$,
and $\Su$ acts on the homology of a weak $\Su$-module~$(M,u)$
by~$\xi_i\cdot[m]=[u(\xi_i)\cdot m]$.
Before describing morphisms of weak modules,
we say how the Koszul functors act on objects.

The Koszul dual of~$(N,v)\in\wMod\Ll$ is defined
as the bounded above $\Su$-module 
$\alg t N=\Su\otimes N$ with differential
\begin{equation}\label{definition-t}
  d(\sigma\otimes n)=\sigma\otimes d n
    +\sum_{\alpha>0}\xi^\alpha\sigma\otimes v(y_\alpha)\cdot n.
\end{equation}
(This is well-defined because $N$ is bounded above.)

The Koszul dual of~$(M,u)\in\wMod\Su$
is the bounded above $\Ll$-module~$\alg h M=\Lu\otimes M$ with differential
\begin{equation}\label{definition-h}
  d(\alpha\otimes m)=(-1)^{|\alpha|}\alpha\otimes d m
    +\sum_{\pi\ne\emptyset}(-1)^{\degree{x_\pi}}
        x_\pi\cdot\alpha\otimes u(\xi_\pi)\cdot m
\end{equation}
and $\Ll$-action coming from that on~$\Lu$,
which is defined similarly to~\eqref{action-cochains},
$$
  (a\cdot\alpha)(a')=(-1)^{\degree a(\degree\alpha+1)}\alpha(a\wedge a').
$$

A morphism~$f$ between two weak $\Ll$-modules $N$~and~$N'$
is a morphism of (strict) $\Su$-modules~$\alg t N\to\alg t N'$.
Its ``base-component''
$$
  N=1\otimes N\hookrightarrow\Su\otimes N
    \stackrel f\longrightarrow\Su\otimes N'
    \to\mkern-14mu\to1\otimes N'=N'
$$
is a chain map inducing
a $\Ll$-equivariant map in homology. If the latter is an isomorphism,
we say that $f$ is a quasi-isomorphism.
The definitions for weak $\Su$-modules are analogous.
The Koszul dual of a morphism of weak modules is what one expects.

The Koszul functors preserve quasi-isomorphisms and are quasi-inverse
to each other,
\cf~\cite[Sec.~2.6]{Franz:03}.
Note that our (left) weak $\Su$-modules correspond
to \emph{left} weak $\Sl$-comodules
and not to right ones as used in~\cite{Franz:03}.
This detail, which is crucial for the present paper, does not affect
Koszul duality.

\medbreak

In the rest of this section we generalise results
of~\cite[Sec.~9]{GoreskyKottwitzMacPherson:98} to weak modules.

\smallskip

A weak $\Su$-module~$M$ is called \emph{split and extended}
if it is quasi-isomorphic to its homology and if the
latter is of the form $\Su\otimes L$ for some graded $R$-module~$L$.
If $M$ it is quasi-isomorphic to its homology
and if the $\Su$-action on~$H(N)$ is trivial,
we say that $M$ is \emph{split and trivial}.
Similar definitions apply to weak $\Ll$-modules.
(Note that it does not make a difference whether we require
the homology of a split and free $\Ll$-module
to be isomorphic to~$\Ll\otimes L$ or to~$\Lu\otimes L$.)

\begin{proposition}\label{duality-trivial-extended}
  Under Koszul duality, split and trivial weak modules
  correspond to split and extended ones.
\end{proposition}

\begin{proof}
  That the Koszul functors carry split and trivial weak modules
  to split and extended ones is almost a tautology.
  The other direction follows from the fact
  that the Koszul functors are quasi-inverse to each other
  and preserve quasi-isomorphisms
  because a split and extended weak module is by definition
  quasi-isomorphic to the Koszul dual of a module with zero differential
  and trivial action.
\end{proof}

\begin{proposition}\label{extended-is-split}
  Let $M$ be in $\wMod\Su$. If $H(M)$ is extended, then $M$ is split.
\end{proposition}

\begin{proof}
  We may assume that $M$ has a strict $\Su$-action because
  any weak $\Su$-module~$M$ is quasi-isomorphic to a strict one
  (for instance, to~$\alg t\alg h M$).
  By assumption, $H(M)\cong\Su\otimes L$ for some graded $R$-module~$L$.
  Since we work over a principal ideal domain, there exists a free resolution
  \[
    0\longleftarrow L\longleftarrow
      P^0\longleftarrow P^1\longleftarrow0
  \]
  of~$L$ with $P^0$,~$P^1$ bounded above.
  Tensoring it with~$\Su$ gives a free resolution of
  the $\Su$-module~$H(M)$
  and therefore the (not uniquely determined) $\Su$-equivariant vertical maps
  in the following commutative diagram with exact rows:
  \begin{diagram*}
    0 & \lTo & \Su\otimes L & \lTo & \Su\otimes P^0 & \lTo & \Su\otimes P^1
      & \lTo 0 \\
    & & \dTo<\cong & & \dTo  & & \dTo \\
    0 & \lTo & H(M) & \lTo & Z(M) & \lTo^d & M .
  \end{diagram*}
  This implies that the 
  total complex~$\Su\otimes P$ is
  quasi-isomorphic to both $H(M)$~and~$M$.
\end{proof}

\section{Proof of Theorem~\protect\ref{A-infinity-formal}}\label{proof-formal}

In this section all algebras are $\N$-graded and connected
unless otherwise stated.

Recall that an $\Ainfty$~map~$f\colon A\shmto A'$
between two algebras is a map of coalgebras~$B A\to B A'$,
see \cite[Sec.~8.1]{McCleary:01}~or~\cite{Keller:01} for example.
It is called strict if it is induced from an algebra map~$A\to A'$.
If $f\colon A\shmto A'$ is $\Ainfty$, then its base component~%
$f_1\colon B_1 A 
\to B_1 A'$ between the elements of external degree~$1$
is a chain map, multiplicative up to homotopy.
We denote the induced algebra map in homology by~$H(f)\colon H(A)\to H(A')$.
If it is an isomorphism, we call $f$ an $\Ainfty$~quasi-isomorphism. 

In order to prove Theorem~\ref{A-infinity-formal},
it is sufficient to construct
an $\Ainfty$~quasi-iso\-mor\-phism~%
$A\shmto H(A)=\bigwedge(x_1,\ldots,x_r)=\Ll$,
due to the following result:

\begin{lemma}\label{shm-inverse}
  Let $A$ be an algebra with $A$~and~$H(A)$ free over~$R$,
  and let $f\colon A\shmto H(A)$ be an $\Ainfty$~map inducing the identity
  in homology.
  Then $f$ has an $\Ainfty$~quasi-inverse, \ie,
  there is an $\Ainfty$~map~$g\colon H(A)\shmto A$ also inducing the identity
  in homology.
\end{lemma}

(At least over fields one can do better:
there any $\Ainfty$~quasi-isomorphism between two algebras
-- even $\Ainfty$~algebras --
is an $\Ainfty$~homotopy equivalence,
\cf~\cite{Munkholm:76} or \cite[Sec.~3.7]{Keller:01}.)

\begin{proof}
  According to~\cite[Prop.~2.2]{Munkholm:74},
  the claim is true if $f$ is strict.
  (Here we use that over a principal ideal domain any quasi-isomorphism
  $A\to H(A)$ of free modules comes from
  a ``trivialised extension'' in the sense of~\cite[\S2.1]{Munkholm:74}.)
  To reduce the general case to this,
  we consider the cobar construction~$\Omega B A$ of~$B A$.
  Coalgebra maps~$h\colon B A\to B A'$ correspond bijectively
  to algebra maps~$\tilde h\colon\Omega B A\to A'$.
  For $h$, the identity of~$A$,
  the map~$\tilde h$ is a quasi-isomorphism
  \cite[Thm~II.4.4]{HusemollerMooreStasheff:74}
  with quasi-inverse (in the category of complexes),
  the canonical inclusion~$A\hookrightarrow\Omega B A$.
  The composition of this map with~$\tilde f\colon\Omega B A\to H(A)$
  is essentially $f_1$, 
  which is a quasi-isomorphism by hypothesis.
  Hence $\tilde f$ is so, too.
  Now compose any $\Ainfty$~quasi-inverse of it
  with the projection~$\Omega B A\to A$.
\end{proof}

Recall that for any complex~$C$ a cycle in~$C^q=\Hom_{-q}(C,R)$
is the same as a chain map~$C\to R[-q]$.
(Here $R[-q]$ denotes the complex~$R$, shifted to degree~$q$.)
The crucial observation, made in a topological context by
Stasheff~\cite[Thm.~5.1]{Stasheff:63b}, is the following:

\begin{lemma}
  $\Ainfty$~maps~$A\shmto\bigwedge(x)$, $\degree x=q>0$,
  correspond bijectively to cocycles in~$(B A)^{q+1}$.
\end{lemma}


\begin{proof}
  Note that the augmentation ideal of~$\bigwedge(x)$ is~$R[-q]$
  (with vanishing product).
  An $\Ainfty$~map~$f\colon A\shmto\bigwedge(x)$ is given
  by components~$f_p\colon\bar A^{\otimes p}\to R[-q]$ of degree~$p-1$
  such that for all~$[a_1,\ldots,a_p]\in B_p(A)$,
  $$
    f_p(d[a_1,\ldots,a_p])=-f_{p-1}(\delta[a_1,\ldots,a_p]),
  $$
  where $d\colon B_p(A)\to B_p(A)$ denotes
  the ``internal'' differential
  and $\delta\colon B_p(A)\to B_{p-1}(A)$ the ``external'' one,
  \cf~\cite[Thm.~8.18]{McCleary:01}.
  In other words, $d(f_p)=-\delta(f_{p-1})$, where $\delta$~and~$d$ now
  denote the dual differentials.
  But this is the condition for a cycle
  in the double complex~$((B A)^*, d, \delta)$ dual to~$BA$.
\end{proof}

By our assumptions, $H^*(B A)=\Su$ is a (negatively graded)
polynomial algebra. Taking representatives of the generators~$\xi_i$ gives
$\Ainfty$~maps~$f^{(i)}\colon A\shmto\bigwedge(x_i)$.
By~\cite[Prop.~3.3~\&~3.7]{Munkholm:74}, they assemble into an $\Ainfty$~map
$$
  f^{(1)}\otimes\cdots\otimes f^{(r)}\colon
  A^{\otimes r}\shmto\bigwedge(x_1)\otimes\cdots\otimes\bigwedge(x_r)=\Ll
$$
whose base component is the tensor product of the base components~$f^{(i)}_1$.
Since $A$ is a Hopf algebra,
the $r$-fold diagonal~$\Delta^{(r)}\colon A\to A^{\otimes r}$
is a morphism of algebras.
A test on the generators~$x_i$ reveals
that the composition~%
$(f^{(1)}\otimes\cdots\otimes f^{(r)})\Delta^{(r)}\colon A\shmto\Ll$
induces an isomorphism in homology, hence is the $\Ainfty$~quasi-isomorphism
we are looking for.

\begin{remark}\rm
  Since we have not really used the coassociativity of~$\Delta$,
  Theorem~\ref{A-infinity-formal} holds even for quasi-Hopf algebras
  in the sense of~\cite[\S IV.5]{HusemollerMooreStasheff:74}.
\end{remark}

\section{The twisting cochain~$v\colon\Sl\to C(G)$}\label{start-other-proofs}

This is now easy: Compose the map~$\Sl\to B\Ll$
determined by the canonical twisting cochain~$w\colon\Sl\to\Ll$
with the map~$B\Ll\to B C(G)$.
This corresponds to a twisting cochain~$\Sl\to C(G)$
mapping each cogenerator~$y_i\in\Sl$ to a representative of~$x_i\in\Ll$.
Since these elements are used to define the $\Ll$-action in the homology
of a weak $\Ll$-module, we get the usual action of~$\Ll=H(G)$ there.

Note that 
by dualisation we obtain
a quasi-isomorphism of algebras~$(B C(G))^*\to\Su$.
This is not exactly the same as
the quasi-isomorphism of algebras~$C^*(BG)\to\Su$
from Proposition~\ref{multiplicative-quasi-isomorphism},
which we are going to construct next.

\section{Proof of Theorem~\protect\ref{main-result} (first part)
  and of Proposition~\protect\ref{multiplicative-quasi-isomorphism}}

In this section we construct maps
$$
  \tcomp_X\colon\Sl\otimes_v C(X)
    \to C(EG\timesunder G X)=C(\Sim t X),
$$
natural in~$X\in\Space G$.
We will show that $\Bf:=\tcomp_{\oneptspace}\colon\Sl\to C(BG)$
is a quasi-isomorphism of coalgebras and that $\tcomp_X$, which maps
from an $\Sl$-comodule to a $C(BG)$--comodule, is a $\Bf$-equivariant
quasi-isomorphism.
Taking duals then gives Proposition~\ref{multiplicative-quasi-isomorphism}
and the first half of Theorem~\ref{main-result}.

Recall that the differential on~$\Sl\otimes_v C(X)$ is
$$
  d(y_\alpha\otimes c)=y_\alpha\otimes d c
    +\sum_{\beta<\alpha}y_\beta\otimes c_{\alpha-\beta}\cdot c,
$$
where we have abbreviated $v(y_{\alpha-\beta})$ to~$c_{\alpha-\beta}$.
The summation runs over all~$\beta$ strictly smaller than~$\alpha$
in the canonical partial ordering of~$\N^r$.

To begin with, we define a map
$$
  f\colon\Sl\otimes_v C(G)\to C(EG)
$$
by recursively setting
\begin{subequations}
\begin{align*}
  f(1\otimes a) &= e_0\cdot a,\\
  f(y_\alpha\otimes a)
    &= \Bigl(S f\bigl(d(y_\alpha\otimes1)\bigr)\Bigr)\cdot a
\end{align*}
\end{subequations}
for~$\alpha>0$.
Here $e_0$ is the canonical base point of the simplicial construction
of the right $G$-space~$EG$
and $S$ its canonical contracting homotopy, \cf~\cite[Sec.~3.7]{Franz:03}.

\begin{lemma}\label{f-quasi-iso}
  This $f$ is a quasi-morphism of right $C(G)$-modules.
\end{lemma}

\begin{proof}
  The map is equivariant by construction.
  By induction, one has for~$\alpha>0$
  $$
    d\,f(y_\alpha\otimes1)
    = d\,S f\bigl(d(y_\alpha\otimes1)\bigr)
    = f\bigl(d(y_\alpha\otimes1)\bigr)-S\,d\,f\bigl(d(y_\alpha\otimes1)\bigr)
    = f\bigl(d(y_\alpha\otimes1)\bigr),
  $$
  which shows that it is a chain map.
  That it induces an isomorphism in homology follows
  from the acyclicity of~$\Sl\otimes_v C(G)$:
  Filter the complex according to the number of factors~$\xi_i$
  appearing in an element~$\xi^\alpha\otimes a$, \ie,
  by~$\alpha_1+\cdots+\alpha_r$.
  Then the $E^1$~term of the corresponding spectral sequence
  is the Koszul complex~$\Sl\otimes_w\Ll$, hence acyclic.
\end{proof}


We will also need the following result:

\begin{lemma}\label{relation-AW-f}
  The image of~$f(y_\alpha\otimes1)$, $\alpha\in\N^r$,
  under the diagonal~$\Delta$ of the coalgebra~$C(EG)$ is
  $$
    \Delta f(y_\alpha\otimes1)
    \equiv\sum_{\beta+\gamma=\alpha}
      f(y_\beta\otimes1)\otimes f(y_\gamma\otimes1),
  $$
  up to terms of the form~$c\cdot a\otimes c'$
  with~$c$,~$c'\in C(EG)$ and~$a\in C(G)$, $\degree a>0$.
\end{lemma}

\begin{proof}
  We proceed by induction, the case~$\alpha=0$ being trivial.
  For~$\alpha>0$ we have
  \begin{align*}
    \Delta f(y_\alpha\otimes1)
    &= \Delta S f\bigl(d\,(y_\alpha\otimes1)\bigr)
     = \sum_{\beta<\alpha}
         \Delta S\bigl(f(y_\beta\otimes1)\cdot c_{\alpha-\beta}\bigr) \\
  \intertext{We now use
    the identity~$\Delta S(c)=S c\otimes 1+(1\otimes S)\AW(c)$
    \cite[Prop.~3.8]{Franz:03}
    and the $C(G)$-equivariance of the Alexander--Whitney map to get}
    &= f(y_\alpha\otimes1)\otimes1
       +(1\otimes S)\sum_{\beta<\alpha}
         \Delta f(y_\beta\otimes1)\cdot\Delta c_{\alpha-\beta}, \\
  \intertext{where the second diagonal is of course that of~$C(G)$.
    By induction and the fact that
    $\Delta c_{\alpha-\beta}\equiv1\otimes c_{\alpha-\beta}$
    up to terms~$a\otimes a'$ with~$\degree a>0$, we find}
    &\equiv f(y_\alpha\otimes1)\otimes1
       +(1\otimes S)\sumsss{\beta+\gamma<\alpha}
         f(y_\beta\otimes1)\otimes
           f(y_\gamma\otimes1)\cdot c_{\alpha-(\beta+\gamma)} \\
    &= f(y_\alpha\otimes1)\otimes1
       +\sumsss{\beta<\alpha \\ \gamma<\alpha-\beta}
         f(y_\beta\otimes1)\otimes
           S f(y_\gamma\otimes c_{(\alpha-\beta)-\gamma}), \\
    \intertext{which simplifies by the definition of~$f$ to}
    &= f(y_\alpha\otimes1)\otimes1+\sum_{\beta<\alpha}
         f(y_\beta\otimes1)\otimes f(y_{\alpha-\beta}\otimes1),
  \end{align*}
  as was to be shown.
\end{proof}

For a $G$-space~$X$ we define the map
$$
  \tcomp_X\colon
    \alg t C(X)=\Sl\otimes_v C(X)\to C(EG\timesunder G X)=C(\Sim t X)
$$
as the bottom row of the commutative diagram
\begin{diagram*}
  & & \Sl\otimes_v C(G)\otimes C(X) & \rTo^{f\otimes1} & C(EG)\otimes C(X)
    & \rTo^\shuffle & C(EG\times X) \\
  & & \dTo & & \dTo & & \dTo \\
  \Sl\otimes_v C(X) & \rEqual
    & \Sl\otimes_v C(G)\otimesunder{C(G)}C(X) & \rTo
    & C(EG)\otimesunder{C(G)}C(X) & \rTo
    & C(EG\timesunder G X),
\end{diagram*}
where $\shuffle$ denotes the shuffle map.
$\tcomp_X$ is obviously natural in~$X$.

It follows from the preceding lemma 
that $\Bf=\tcomp_{\oneptspace}\colon\Sl\to C(BG)$ is a morphism of coalgebras
because terms of the form~$c\cdot a$ with~$\degree a>0$ are annihilated by
the projection~$C(EG)\to C(BG)$. (We are working with normalised chains!)
Using naturality and the commutativity of the diagram
\begin{diagram*}
  C(EG)\otimes C(X) & \rTo^\shuffle & C(EG\times X) \\
  \dTo<{\Delta_{C(EG)}\otimes1} & & \dTo>{\Delta_{C(EG\times X)}} \\
  C(BG)\otimes C(EG)\otimes C(X) & \rTo_{1\otimes\shuffle}
    & C(BG)\otimes C(EG\times X),
\end{diagram*}
one proves similarly
that $\tcomp_X$ is a $\Bf$-equivariant morphism of comodules.
To see that it induces an isomorphism in homology,
consider the diagram
\begin{diagram*}
  \Tor^{C(G)}\bigl(\Sl\otimes_v C(G),C(X)\bigr)
    & \rTo & H\bigl(\Sl\otimes_v C(G)\otimesunder{C(G)}C(X)\bigr)
    & \rEqual & H(\Sl\otimes_v C(X))\\
  \dTo<{\Tor^{\id}(f,\id)} & & \dTo & & \dTo>{H(\tcomp_X)} \\
  \Tor^{C(G)}\bigl(C(EG), C(X)\bigr)
    & \rTo & H\bigl(C(EG)\otimesunder{C(G)}C(X)\bigr)
    & \rTo & H(EG\timesunder G X).
\end{diagram*}
The composition along the bottom row is an isomorphism
by Moore's theorem~\cite[Thm.~7.27]{McCleary:01},\footnote{In fact,
each single arrow is an isomorphism.
This follows from the twisted Eilenberg--Zilber theorem,
see \cite{Gugenheim:72} for example.}
and the top row is so because $\Sl\otimes_v C(G)$ is $C(G)$-flat.
Since $\Tor^{\id}(f,\id)$ is an isomorphism by Lemma~\ref{f-quasi-iso},
$H(\tcomp_X)$ is so, too.

\section{The twisting cochain~$u\colon\Lu\to C^*(BG)$ \\
  and the end of the proof of
  Theorem~\protect\ref{main-result}}\label{end-other-proofs}

The map~$\Bf\colon\Sl\to C(BG)$ is a quasi-isomorphism
of simply connected coalgebras.
Similar to the first step in the proof of Lemma~\ref{shm-inverse},
it comes from a trivialised extension (or ``Eilenberg--Zilber data''
in the terminology of~\cite{GugenheimMunkholm:74}).
By~\cite[Thm.~4.1$^*$]{GugenheimMunkholm:74},
there is an algebra map~$F\colon\Omega C(BG)\to\Omega\Sl$
whose base component~$F_{-1}\colon\Omega_{-1}C(BG)\to\Omega_{-1}\Sl$
is essentially the chosen homotopy inverse to~$\Bf$.
Composing such an~$F$ with the canonical map~$g\colon\Omega\Sl\to\Ll$,
we get a twisting cochain~$\tilde u\colon C(BG)\to\Ll$.
Write
\begin{equation}\label{cochain-u-tilde}
  \tilde u=\sumsss{\emptyset\ne\pi\subset\{1,\ldots,r\}}x_\pi\otimes\gamma_\pi
  \in\Ll\otimes C^*(BG)=\Hom(C(BG),\Ll).
\end{equation}
Then $\gamma_i$ is a representative of the generator~$\xi_i\in\Su$
because it is a cocycle (\cf~\cite[eq.~(2.12)]{Franz:03}) and
$$
  \tilde u(\Bf(y_i))=g(F([\Bf(y_i)]))=g(y_i)=x_i.
$$

The dual~$u=\tilde u^*\colon\Lu\to C^*(BG)$ is again a cochain,
which corresponds under the isomorphism~$\Hom(\Lu,C^*(BG))=C^*(BG)\otimes\Ll$
to the transposition of factors of~\eqref{cochain-u-tilde}.
Therefore, the induced action of~$\Su$ on a $C^*(BG)$-module, considered
as weak $\Su$-module, is given by~$\xi_i\cdot[m]=[\gamma_i\cdot m]$,
as desired.

\medbreak

For a given $G$-space~$X$, we now look at the map~$\tcomp_X^*$
as a quasi-isomorphism of $C^*(BG)$-modules,
where the module structure of~$\alg t C^*(X)$
is induced by~$\Bf^*$.
By naturality, it is a morphism of weak $\Su$-modules.
This new weak $\Su$-action on~$\alg t C^*(X)$
coincides with the (strict) old one because the composition
$$
  (\Omega\Sl)^*\stackrel{F^*}\longrightarrow(\Omega C(BG))^*
    \stackrel{\Bf^*}\longrightarrow(\Omega\Sl)^*
$$
is the identity.
This proves that $\tcomp_X^*$ is a quasi-isomorphism of weak $\Su$-modules,
hence that the functors $C^*\circ\Sim t$~and~$\alg t\circ C^*$
are quasi-isomorphic.

The corresponding result for the functors $\Sim h$~and~$\alg h$
is a formal consequence of this
because they are quasi-inverse to $\Sim t$~and~$\alg t$, respectively.
This finishes the proof of Theorem~\ref{main-result}.

\medbreak

\begin{remark}\rm
  For~$G=(S^1)^r$ a torus (and a reasonable choice of~$v$)
  one may also take the twisting cochain~$\Lu\to C^*(BG)$
  of Gugenheim--May~\cite[Example~2.2]{GugenheimMay:74},
  which is defined using
  iterated \cuponeproduct s of (any choice of)
  representatives~$\gamma_i\in C^*(BG)$
  of the~$\xi_i\in\Su$.
  (This follows for example from~\cite[Cor.~4.4]{Franz:03}.) 
  It would be interesting to know whether this remains true in general
  if one chooses the~$\gamma_i$ carefully enough.
\end{remark}

\section{Equivariantly formal spaces}\label{equivariant-cohomology}

An important class of $G$-spaces are the equivariantly formal ones.
Their equivariant cohomology is particularly simple, which is often
exploited in algebraic or symplectic geometry or combinatorics.

We say that $X$ is \emph{$R$-equivariantly formal}
if the following conditions hold.

\begin{proposition}
  For a $G$-space~$X$, the following are equivalent:
  \begin{enumerate}
  \item \label{condition-extended}
    $H_G^*(X)$ is extended.
  \item \label{condition-split-and-extended}
    $C^*(X_G)$ is split and extended.
  \item \label{condition-split-and-trivial}
    $C^*(X)$ is split and trivial.
  \item \label{condition-section}
    The canonical map~$H_G^*(X)\to H^*(X)$ admits
    a section of graded $R$-modules.
  \item \label{condition-E2}
    $H_G^*(X)$ is isomorphic, as $\Su$-module,
    to the $E_2$~term~$\Su\otimes H^*(X)$
    of the Leray--Serre spectral sequence for~$X_G$
    (which therefore degenerates).
  \end{enumerate}
\end{proposition}

Note that if $R$ is a field,
condition~\eqref{condition-extended} means that $H_G^*(X)$ is free over~$\Su$,
and condition~\eqref{condition-section} that $H_G^*(X)\to H^*(X)$
is surjective. A space~$X$ with the latter property
is traditionally called ``totally non-homologous to zero in~$X_G$
with respect to~$R$''.
We stress the fact that for some of the above conditions
we really need the assumption that $R$ is a principal ideal domain.

In~\cite{FranzPuppe:??} it is shown that a compact symplectic manifold~$X$
with a Hamiltonian torus action is $\Z$-equivariantly formal
if $X^T=X^{T_p}$ for each prime~$p$ that kills elements in~$H^*(X^T)$.
Here $T_p\cong\Z_p^r$ denotes the maximal $p$-torus contained in the torus~$T$.
In particular, a compact Hamiltonian $T$-manifold is $\Z$-equivariantly formal
if the isotropy group of each non-fixed point is contained
in a proper subtorus.

\begin{proof}
  $\eqref{condition-E2}\Rightarrow\eqref{condition-extended}$ is trivial.
  $\eqref{condition-extended}\Rightarrow\eqref{condition-split-and-extended}$
  follows from Proposition~\ref{extended-is-split}, and
  $\eqref{condition-split-and-extended}\Rightarrow\eqref{condition-split-and-trivial}$
  from Proposition~\ref{duality-trivial-extended}
  because $C^*(X)$ and $C^*(X_G)$ are Koszul dual by Theorem~\ref{main-result}.
  $\eqref{condition-section}\Rightarrow\eqref{condition-E2}$
  is the Leray--Hirsch theorem.
  (Note that it holds here for arbitrary~$X$
  because $H^*(BG)=\Su$ is of finite type.)

  $\eqref{condition-split-and-trivial}\Rightarrow\eqref{condition-section}$:
  The (in the simplicial setting canonical)
  map~$C^*(X_G)\to C^*(X)$ is the composition of~$\tcomp_X^*$
  with the canonical projection~$\alg t C^*(X)\to C^*(X)$.
  Since $C^*(X)$ is split, we can pass from~$C^*(X)$ to~$H^*(X)$
  by a sequence of commutative diagrams
  \begin{diagram*}
    \alg t N & \rTo & N \\
    \dTo & & \dTo \\
    \alg t N' & \rTo & N' \\
  \end{diagram*}
  where the vertical arrow on the right is the base component
  of the quasi-isomorphism of weak $\Ll$-modules given on the left.
  But for the projection~$\Su\otimes H^*(X)\to H^*(X)$ the assertion
  is obvious because $\Ll$ acts trivially on~$H^*(X)$, which means
  that there are no differentials any more.
\end{proof}

\section{Relation to the Cartan model}\label{other-models}



In differential geometry and differential homological algebra
many different complexes (``models'') are known that compute
the equivariant cohomology of a space.
We content ourselves with indicating the relation between our construction
and the probably best-known one, the so-called Cartan model.
We use real or complex coefficients.

Let $G$ be a compact connected Lie group and $X$ a $G$-manifold.
The Cartan model of~$X$ is the complex
\begin{subequations}\label{definition-Cartan-model}
\begin{equation}
  \Bigl(\Sym(\g^*)\otimes\Omega(X)\Bigr)^G  
\end{equation}
of $G$-invariants with differential
\begin{equation}\label{differential-Cartan-model}
  d(\sigma\otimes\omega)
  =\sigma\otimes d\omega+\sum_{j=1}^s\zeta_j\sigma\otimes z_j\cdot\omega.
\end{equation}
\end{subequations}
Here $\Sym(\g^*)$ denotes the (evenly graded) polynomial functions
on the Lie algebra~$\g$ of~$G$,
$(z_j)$ a basis of~$\g$ with dual basis~$(\zeta_j)$,
and $z_j\cdot\omega$ the contraction of the form~$\omega$
with the generating vector field associated with~$z_j$.
The Cartan model computes $H_G^*(X)$ as algebra and as $\Su$-module,
\cf~\cite{GuilleminSternberg:99}.

As alluded to in the introduction,
Goresky, Kottwitz and MacPherson~\cite{GoreskyKottwitzMacPherson:98}
have found an even smaller complex giving the $\Su$-module~$H_G^*(X)$,
namely $\alg t\Omega(X)^G$, or explicitly
\begin{subequations}\label{definition-small-model}
\begin{equation} 
  \Sym(\g^*)^G\otimes \Omega(X)^G,
\end{equation}
where $\Omega(X)^G$ denotes the complex of $G$-invariant differential forms
on~$X$. The differential
\begin{equation}
  d(\sigma\otimes\omega)
  =\sigma\otimes d\omega+\sum_{i=1}^r\xi_i\sigma\otimes x_i\cdot\omega
\end{equation}
\end{subequations}
is similar to~\eqref{differential-Cartan-model},
but the summation now runs over a system of generators
of~$\Su=H^*(BG)=\Sym(\g^*)^G$.
(This is of course differential~\eqref{definition-t} for strict $\Ll$-modules.)
Maszczyk and Weber~\cite{MaszczykWeber:02} have proved that the
complexes \eqref{definition-Cartan-model}~and~\eqref{definition-small-model}
are quasi-isomorphic as $\Su$-modules.

For the case of torus actions
(where \eqref{definition-Cartan-model}~and~\eqref{definition-small-model}
coincide), Goresky--Kottwitz--MacPherson~%
\cite[Sec.~12]{GoreskyKottwitzMacPherson:98}
have shown that one may replace $\Omega(X)^T$
by singular cochains together with the ``sweep action'',
which is defined by restricting the action of~$C(T)$ along
a quasi-isomorphism of algebras~$\Ll=H(T)\to C(T)$.
The latter is easy to construct: Since $C(T)$ is (graded) commutative,
any choice of representatives~$c_i$ of a basis~$(x_i)$ of~$H_1(T)$
can be extended to such a morphism.
Now all ingredients are defined for an arbitrary topological $T$-space~$X$
and an arbitrary coefficient ring~$R$, and the resulting complex does
indeed compute $H_T^*(X)$ as algebra and as $\Su$-module
in this generality,
see F\'elix--Halperin--Thomas \cite[Sec.~7.3]{FelixHalperinThomas:95}.

\section*{Appendix: The theorems of Samelson--Leray and Hopf}

All differentials are zero in this section.
Recall that an element~$a$ of a Hopf algebra~$A$ is called primitive
if $\Delta a=a\otimes 1+1\otimes a$ or, equivalently,
if the projection of~$\Delta a$ to~$\bar A\otimes\bar A$ is zero.

Let $A$ be a Hopf algebra over a field, isomorphic as algebra
to an exterior algebra. Then $A$ is primitively generated (Samelson--Leray).
If $R$ is a field of characteristic~$0$ and $A$ a connected
commutative Hopf algebra, finite-dimensional over~$R$,
then multiplicatively it is an exterior algebra (Hopf),
hence also primitively generated.
(A good reference for our purposes is \cite[\S\S1,~2]{Kane:88}.)

We now show that the analogous statements
hold over any principal ideal domain.
Denote for a Hopf algebra~$A$ over~$R$
the extension of coefficients to the quotient field of~$R$ by~$A_{(0)}$.

\begin{proposition}
  Let $A$ be a Hopf algebra, free over~$R$
  and such that $A_{(0)}$ is a primitively generated exterior algebra.
  Then $A$ is a primitively generated exterior algebra.
\end{proposition}

\begin{proof}
  Let $A'$ be the sub Hopf algebra generated by the free submodule
  of primitive elements of~$A$. Then $A'_{(0)}=A_{(0)}$ (Samelson--Leray),
  hence $A'$ is a primitively generated exterior algebra
  and $A/A'$ is $R$-torsion.
  Take an~$a\in A\setminus A'$ of smallest degree.
  Then $k a\in A'$ for some~$0\ne k\in R$, and the image of~$\Delta a$
  in~$\bar A\otimes\bar A$ already lies in~$\bar A'\otimes\bar A'$.
  Write $k a=a_1+a_2$ with~$a_1\in A'$ primitive
  and~$a_2\in\bar A'\cdot\bar A'$.
  Note that the image of~$\Delta a_2$ in~$\bar A'\otimes\bar A'$
  is divisible by~$k$.
  This implies that $a_2$ is divisible by~$k$ in~$A'$.
  (Look at how the various products of the generators
  of a primitively generated exterior algebra behave under the diagonal.)
  Since $a-a_2/k$ is primitive, it lies in~$A'$, hence $a$ as well.
  Therefore, $A=A'$.
\end{proof}

\bibliographystyle{plain}
\bibliography{abbrev,algtop}

\begin{thebibliography}{10}

\bibitem{BernsteinGelfandGelfand:78}
I.~N. Bern{\v{s}}te{\u\i}n, I.~M. Gel'fand, and S.~I. Gel'fand.
\newblock Algebraic vector bundles on~{${\mathbf P}^n$} and problems of linear
  algebra.
\newblock {\em Funkt. Anal. Appl.} {\bf 12} (1978), 212--214.


\bibitem{DrorDwyerKan:80}
E. Dror, W.~G. Dwyer, D.~M. Kan.
\newblock Equivariant maps which are self homotopy equivalences.
\newblock {\em Proc. AMS} {\bf 80} (1980), 670--672.

\bibitem{FelixHalperinThomas:95}
Y.~Felix, S.~Halperin, and J.-C.~Thomas.
\newblock Differential graded algebras in topology.
\newblock Ch.~16 in: I.~M.~James (ed.), {\em Handbook of algebraic topology},
  Elsevier, Amsterdam 1995.

\bibitem{Franz:03}
M.~Franz.
\newblock {K}oszul duality and equivariant cohomology for tori.
\newblock {\em Int.\ Math.\ Res.\ Not.} {\bf 42}, 2255-2303 (2003).

\bibitem{FranzPuppe:??}
M.~Franz and V.~Puppe.
\newblock Exact sequences for equivariantly formal spaces.
\newblock ArXiv e-print math.AT/0307112.

\bibitem{GoreskyKottwitzMacPherson:98}
M.~Goresky, R.~Kottwitz, and R.~MacPherson.
\newblock Equivariant cohomology, {K}oszul duality, and the localization
  theorem.
\newblock {\em Invent.\ Math.} {\bf 131} (1998), 25--83.

\bibitem{Gugenheim:72}
V.~K.~A.~M.~Gugenheim.
\newblock On the chain-complex of a fibration. 
\newblock {\em Illinois J.\ Math.} {\bf 16} (1972), 398--414.
 

\bibitem{GugenheimMay:74}
V.~K.~A.~M. Gugenheim and J.~P.~May.
\newblock On the theory and applications of differential torsion products.
\newblock {\em Mem.\ Am.\ Math.\ Soc.} {\bf 142} (1974).

\bibitem{GugenheimMunkholm:74}
V.~K.~A.~M. Gugenheim and H.~J. Munkholm.
\newblock On the extended functoriality of $\rm Tor$~and~$\rm Cotor$.
\newblock {\em J.\ Pure Appl.\ Algebra} {\bf 4} (1974), 9--29.

\bibitem{GuilleminSternberg:99}
V.~W. Guillemin and S.~Sternberg.
\newblock {\em Supersymmetry and equivariant {d}e {R}ham theory}.
\newblock 
  Springer, Berlin 1999.

\bibitem{Huebschmann:84}
J.~Huebschmann.
\newblock Perturbation theory and small models for the chains
  of certain induced fibre spaces.
\newblock Habilitationsschrift, Universit\"at Heidelberg 1984.

\bibitem{HusemollerMooreStasheff:74}
D.~H. Husemoller, J.~C. Moore, and J. Stasheff.
\newblock Differential homological algebra and homogeneous spaces.
\newblock {\em J.\ Pure Appl.\ Algebra} {\bf 5} (1974), 113--185.

\bibitem{Kane:88}
R.~M.~Kane.
\newblock {\em The homology of {H}opf spaces}.
\newblock North-Holland, Amsterdam 1988.

\bibitem{Keller:01}
B.~Keller.
\newblock Introduction to {$A$}-infinity algebras and modules.
\newblock {\em Homology Homotopy Appl.} {\bf 3} (2001), 1--35.
\newblock Addendum {\em ibid.}~{\bf 4} (2002), 25--28.

\bibitem{MaszczykWeber:02}
T.~Maszczyk and A.~Weber.
\newblock {K}oszul duality for modules over {L}ie algebras.
\newblock {\em Duke Math.\ J.} {\bf 112} (2002), 511--520.


\bibitem{McCleary:01}
J.~Mc{C}leary.
\newblock {\em A user's guide to spectral sequences}, 2nd ed.,
\newblock Cambridge University Press, Cambridge 2001.


\bibitem{Munkholm:74}
H.~J. Munkholm.
\newblock The {E}ilenberg--{M}oore spectral sequence and strongly homotopy
  multiplicative maps.
\newblock {\em J.\ Pure Appl.\ Algebra} {\bf 5} (1974), 1--50.

\bibitem{Munkholm:76}
H.~J. Munkholm.
\newblock {shm} maps of differential algebras, {I}. A characterization up to homotopy.
\newblock {\em J.\ Pure Appl.\ Algebra} {\bf 9} (1976), 39--46.

\bibitem{Stasheff:63b}
J.~D.~Stasheff.
\newblock Homotopy associativity of {$H$}-spaces. {II}.
\newblock {\em Trans.\ AMS} {\bf 108} (1963), 293--312.

\bibitem{StasheffHalperin:70}
J.~Stasheff and S.~Halperin.
\newblock Differential algebra in its own rite.
\newblock Pages 567--577
  in {\em Proc. Adv. Study Inst. Alg. Top.
    (Aarhus 1970)}, Vol. III.
\newblock Various Publ. Ser. {\bf 13},
  Mat. Inst. Aarhus Univ., Aarhus 1970.

\end{thebibliography}

\end{document}